\documentclass[12pt]{article}
\usepackage{graphicx}

\usepackage{epsfig}
\usepackage{amsmath}
\usepackage{amssymb}
\usepackage{amsthm}
\usepackage{latexsym}

\newcommand{\comment}[1]{}
\newtheorem{theorem}{Theorem}
\newtheorem{lemma}{Lemma}[section]
\newtheorem{remark}{Remark}[section]
\newtheorem{corollary}{Corollary}[section]
\newtheorem{proposition}{Proposition}[section]

\newtheorem{definition}{Definition}[section]
\begin{document}

\title{{\LARGE\sf
{\bf Scaling Limit and Critical Exponents for Two-Dimensional Bootstrap Percolation}}
}

\author{
{\bf Federico Camia}
\thanks{Research partially supported by a Marie Curie Intra-European Fellowship
under contract MEIF-CT-2003-500740.}\,
\thanks{E-mail: camia\,@\,eurandom.tue.nl}\\
{\small \sl EURANDOM, P.O. Box 513, 5600 MB Eindhoven, The Netherlands}
}

\date{}

\maketitle

\begin{abstract}
Consider a cellular automaton with state space $\{ 0,1 \}^{{\mathbb Z}^2}$
where the initial configuration $\omega_0$ is chosen according to a Bernoulli
product measure, $1$'s are stable, and $0$'s become $1$'s if they are
surrounded by at least three neighboring $1$'s.
In this paper we show that the configuration $\omega_n$ at time $n$
converges exponentially fast to a final configuration $\bar\omega$,
and that the limiting measure corresponding to $\bar\omega$ is in
the universality class of Bernoulli (independent) percolation.

More precisely, assuming the existence of the critical exponents
$\beta$, $\eta$, $\nu$ and $\gamma$, and of the continuum scaling
limit of crossing probabilities for independent site percolation
on the close-packed version of ${\mathbb Z}^2$ (i.e., for independent
$*$-percolation on ${\mathbb Z}^2$), we prove that the
bootstrapped percolation model has the same scaling limit
and critical exponents.

This type of bootstrap percolation can be seen as a paradigm
for a class of cellular automata whose evolution is given,
at each time step, by a monotonic and nonessential
enhancement~\cite{ag,grimmett}.
\end{abstract}

\noindent {\bf Keywords:} bootstrap percolation, scaling limit,
critical exponents, universality.

\noindent {\bf AMS 2000 Subject Classification:} 82B27, 60K35, 82C22,
82B43, 82C20, 82C43, 37B15.

\section{Introduction and Motivations}

Bootstrap percolation is a cellular automaton with
state space $\{ 0, 1 \}^{{\mathbb Z}^d}$ which evolves in
discrete time according to the following rule: a given
configuration of $0$'s and $1$'s is updated by changing
to $1$ each $0$ with at least $l$ neighboring $1$'s and
leaving the rest of the configuration unchanged.
Here $l$ is a nonnegative integer no bigger than $2d$,
and the initial configuration is chosen according to a
Bernoulli product measure with parameter $p$ (the initial
density of $1$'s).

It is known from the work of van Enter~\cite{vanenter}
and Schonmann~\cite{schonmann} that if $l \leq d$, then
almost all initial configurations evolve toward the
constant configuration with $1$'s at all sites.
On the other hand, it is clear that when $l > d$,
the $1$'s do not take over completely, for almost every
initial configuration.
For example, if $l = 2d$, the only $0$'s that become $1$'s
are those completely surrounded by $1$'s.
The configuration changes only once and the final measure
is in some sense very close to a product measure.

The case $l = 2d-1$ is already much more interesting;
it is studied in~\cite{fis}, where the limiting measure
(whose existence is ensured by the monotonicity of the dynamics)
is shown to have exponentially decaying correlations, and
the density function to be analytic in $[0, 1]$
(for simplicity, the authors restrict their attention
to $d=2$, but all arguments used are immediately seen to
hold qualitatively for any $d \geq 2$).
In this paper, we consider the same model studied in~\cite{fis},
with $l = 2d-1$ and $d=2$.

The exponential decay of correlations proved in~\cite{fis}
shows that the bootstrap dynamics generates only short-range
correlations between different sites.
It is an open question, in general, whether introducing
short-range correlations modifies the critical exponents
and the continuum scaling limit (see Section~\ref{universality}).
Based on very general renormalization group arguments,
the answer to this question is expected to be negative
under a broad class of conditions (see, for example, \cite{cardy2}),
but very few rigorous results are available, especially below
the \emph{upper critical dimension}, where the values of the
critical exponents are expected to be different from those
predicted by mean-field theory (there are, however, some
exceptions -- see, e.g., \cite{ps,cns1,cns2,cn1,camia}).
The main goal of this paper is to present a model
for which this question can be answered rigorously.

Our first result, Theorem~\ref{convergence} of
Section~\ref{bootstrap-results}, states that the probability
$\Pi(n)$ that the origin changes state after time $n$ decays
exponentially in $n$.
Following its proof, we present a proof of the exponential
decay of correlations, Theorem~\ref{exp-corr}, which is somewhat
different than that of~\cite{fis}.
The purpose of our proof is to show that the same mechanism
is responsible for the exponential decay of correlations and
the exponential convergence to the final configuration
(Theorem~\ref{convergence}).

The mechanism we are referring to has to do with the nature
of the bootstrap dynamics (which only removes ``dangling ends"
from clusters of $0$'s, starting from the tip) combined with
the fact that large ``tree-like" clusters are ``unlikely" in
Bernoulli percolation.
These two observations are sufficient to prove the results
mentioned above, and also imply that the bootstrap dynamics
only removes relatively small pieces of clusters of $0$'s,
leaving unchanged the large pieces that are relevant in the
scaling limit.
This is the reason why the bootstrapped measure has the same
critical exponents and scaling limit as independent percolation.
In this (very special) case, therefore, the reason for the
``universal" behavior (see Section~\ref{universality}) of the
bootstrap model is quite clear, without having to resort to
more general (non-rigorous) renormalization group arguments.

The main results are presented in Section~\ref{main-results},
and the proofs given in Section~\ref{proofs}.
The proofs are based on ideas developed in~\cite{cns1}
and~\cite{cn1}, and the main tool is the natural coupling
between independent and bootstrap percolation provided
by the bootstrap dynamics itself, which allows to draw
conclusions regarding the bootstrap model by estimating
the probability of events under the initial product measure.
In particular, the coupling allows to compare bootstrap
and Bernoulli percolation and show that they do not differ
``macroscopically."



Besides being interesting in its own right, the bootstrap
dynamics considered in this paper can also be seen as a
particular example of a special class of cellular automata
on various lattices whose evolution is given, at each time
step, by a \emph{monotonic} and \emph{nonessential enhancement}
of \emph{finite range} (see~\cite{ag,grimmett} for the
relevant definitions).
Such cellular automata would be called \emph{subcritical}
in the language of~\cite{gg} (although they do not represent
all subcritical cellular automata).

In order to extend the results of the present paper to the
whole class of cellular automata specified above, one needs
first of all to find candidates for the protected sites
of Definition~\ref{stable+protected} below.
The existence of suitable candidates for that role is not
obvious in that generality, but can be proved using results
of~\cite{camia}.
Once this is done, the proofs of the main theorems would
proceed in much the same way as in this paper.

%
%
%
%

\subsection{Universality} \label{universality}

Our main motivation for studying the type of questions addressed
here (see also~\cite{camia,cn1,cns1,cns2}) is related to the
idea of universality, according to which most statistical-mechanical
systems fall into universality classes such that systems belonging
to the same class have the same critical exponents (the exponents
describing the nature of the divergence of certain quantities or
their derivatives near or at the critical point, where a second
order phase transition occurs).

A closely related notion of universality has to do with the
continuum scaling limit, a limit in which the \emph{microscopic}
scale of the system (e.g., the lattice spacing for systems defined
on a lattice) is sent to zero, while focus is kept on features
manifested on a \emph{macroscopic} scale.
Such a limit is only meaningful at the critical point, where
the correlation length (i.e., the ``natural length scale'' of
the system) is supposed to diverge.
It seems that universality in terms of the scaling limit is
a stronger notion than that in terms of critical exponents.
In~\cite{sw}, some knowledge of the scaling limit is used to
determine critical exponents in the case of two-dimensional
independent site percolation on the triangular lattice, but
there is no general result in that direction.

The concept of universality and the existence of universality
classes arise naturally in the theory of critical phenomena based
on the renormalization group, and are backed by strong theoretical
and experimental evidence.
Below the upper critical dimension, however, only few rigorous
results are available.

\section{Definition of the Model and Preliminary Results} \label{bootstrap-results}

Consider a bootstrap percolation model on ${\mathbb Z}^2$
with initial configuration $\omega = \{ \omega(x) \}_{{x \in \mathbb Z}^2}
\in \{ 0, 1 \}^{{\mathbb Z}^2}$ chosen according to a product
$P_p = \Pi_{x \in {\mathbb Z}^2} \, \nu_x$ of Bernoulli measures
$\{ \nu_x \}_{x \in {\mathbb Z}^2}$ with parameter $p$
(i.e., $\nu_x[\omega(x)=1]=p=1-\nu_x[\omega(x)=0]$); $E_p$ will
denote expectation with respect to $P_p$.
The evolution is given by the following rules:
\begin{itemize}
\item updates are performed at discrete times $n =1, 2, \ldots$
\item $1$'s are stable
\item at the next update, a $0$ becomes $1$ if it has at least three
neighboring $1$'s.
\end{itemize}

Given an initial configuration $\omega$, the bootstrapped configuration
is denoted by $\bar\omega$ and the limiting distribution by $\bar P_p$.
We will call the sites of ${\mathbb Z}^2$ \emph{open} if they are assigned
value $1$ and \emph{closed} if they are assigned value $0$.
Given a subset $D$ of ${\mathbb Z}^2$, we denote $|D|$ its cardinality
and by $\omega_D$ the configuration $\omega$ restricted to $D$.
A subset $D \in {\mathbb Z}^2$ is called a \emph{plaquette} if it is
composed of four sites which are the vertices of a square of side
length $1$.

We denote by $p_c$ the critical value of independent site percolation
on ${\mathbb Z}^2$ and by $p_c^*=1-p_c$ the  critical value of
independent $*$-percolation on the same lattice, which corresponds to
site percolation  on ${\mathbb Z}^2_{cp}$, the \emph{close-packed} version
of ${\mathbb Z}^2$ (obtained by adding the diagonals to each face of
${\mathbb Z}^2$).
We call \emph{${\mathbb Z}^2$-path} (respectively, \emph{$*$-path}) an
ordered sequence $(x_0, \ldots, x_k)$ of sites of ${\mathbb Z}^2$ such
that $x_{i-1}$ and $x_i$ are neighbors in ${\mathbb Z}^2$ (resp., in
${\mathbb Z}^2_{cp}$) for $i = 1, \ldots, k$ and $x_i \neq x_j$
for $i \neq j$.
A \emph{${\mathbb Z}^2$-loop} (resp., \emph{$*$-loop}) is a
${\mathbb Z}^2$-path (resp., $*$-path) that ends at a
${\mathbb Z}^2$-neighbor (resp., $*$-neighbor) of the
starting site.
A path or a loop will be called closed or open if all its sites
are closed or open, respectively.
We call \emph{length} of a path or loop the number of sites in it.

\begin{definition} \label{stable+protected}
A closed site $x \in {\mathbb Z}^2$ is called \emph{stable} if
and only if $\bar\omega(x)=0$.
A site is said to be \emph{protected} if it is closed and is
part of a group of four closed sites forming a plaquette.
\end{definition}

\noindent Clearly, a protected site is stable, together with the other
three sites that complete the plaquette of Definition~\ref{stable+protected},
since each one of them has (at least) two closed ${\mathbb Z}^2$-neighbors.

The following are two elementary but useful lemmas.

\begin{lemma} \label{lemma-protected}
If $x$ and $y$ are stable closed sites and $\omega$ contains a closed
${\mathbb Z}^2$-path $\pi$ joining $x$ and $y$, then all the sites
in $\pi$ are stable.
Closed ${\mathbb Z}^2$-loops are also stable.
\end{lemma}

\noindent {\it Proof.} For the first claim, it is enough to observe that
each site in $\pi$ other than $x$ or $y$ has at least two closed
${\mathbb Z}^2$-neighbors in $\omega$.
In a ${\mathbb Z}^2$-loop, every site has at least two closed
${\mathbb Z}^2$-neighbors. \fbox{} \\

For $(x,x')$ an ordered pair of neighbors in ${\mathbb Z}^2$,
we define the \emph{partial cluster} $C_{(x,x')}$
to be the set of sites $y \in {\mathbb Z}^2$ such that
there is a ${\mathbb Z}^2$-path $(x_0=x', x_1, \dots, x_k=y)$,
with $x_1 \neq x$, whose sites are all open or all closed.

\begin{lemma} \label{stable-path}
A closed ${\mathbb Z}^{2}$-path $(y_0, \dots, y_k)$
in $\omega$ is stable (i.e., all its sites are stable)
if $C_{(y_1,y_0)}$ and $C_{(y_{k-1},y_k)}$ both contain
protected sites.
\end{lemma}

\noindent \emph{Proof.}
The path $(y_0, \dots, y_k)$ in $\omega$ is stable
because there exists a (generally longer) closed
path that starts and ends at stable sites and
contains $(y_0, \dots, y_k)$ as a subpath.
Since the starting and ending sites of such a
path are stable, all the other sites of the path,
including $y_0, \dots, y_k$, are also stable by
an application of Lemma~\ref{lemma-protected}. \fbox{} \\

We will denote by $\omega_n$ the percolation configuration at time
$n$, i.e., after $n$ updates of the initial configuration.
With this notation we have $\omega_0 = \omega$ (the initial configuration)
and $\omega_{\infty} = \bar \omega$ (the final configuration).
Our first result concerns the speed of convergence of $\omega_n$
to $\bar\omega$.

\begin{theorem} \label{convergence}
Let $\Pi(n)$ be the probability that the origin changes state after
time $n$. Then, for each $p \in [0,1]$ there exists $c_0>0$ such that
$\Pi(n) \leq \exp{(-c_0 \, n)}$.
\end{theorem}

\noindent {\it Proof}. Let $o$ denote the origin of ${\mathbb Z}^2$.
If $\omega_0(o)=1$, the origin never changes state, therefore we
will assume, without loss of generality, that $\omega_0(o)=0$ and
also that $0<p<1$.
To analyze when the origin becomes $1$, we consider its cluster
$C_o$ at time $0$.
Let $x_1, x_2, x_3, x_4$ be the four ${\mathbb Z}^2$-neighbors of
the origin in some deterministic order.
For $x_i, \, i=1,2,3,4$, we define the \emph{branch} $C_i$
to be the partial cluster $C_{(o,x_i)}$.
If $\omega_0(x_i)=1$, we say that $C_i$ is \emph{empty}.

Our first observation is that if the branches $C_i, \, i=1,2,3,4$,
are not distinct, the origin belongs to a ${\mathbb Z}^2$-loop and
is stable by Lemma~\ref{lemma-protected}.
We also notice that, for the origin to become $1$, no more than one
branch $C_i$ can have a stable site, otherwise the origin would
again be stable by Lemma~\ref{lemma-protected}.
We will then assume that the branches $C_i, \, i=1,2,3,4$, are
distinct, and that at most one of them contains a stable site.
Notice that the branches that do not contain stable sites have
a tree-like structure (they do not contain ${\mathbb Z}^2$-loops).

Consider first the case in which exactly one branch contains
a stable site.
The origin will then become $1$ at some time $n$ equal to one
plus the length of a longest self-avoiding ${\mathbb Z}^2$-path
contained in one of the remaining branches.
If no branch contains a stable site, let $C_j$ be a branch
containing a longest ${\mathbb Z}^2$-path and $\pi$ be a longest
${\mathbb Z}^2$-path not contained in $C_j$.
Then the origin will become $1$ at some time $n$ equal to one
plus the length of $\pi$.

The discussion above shows that a necessary condition for the
origin to change state after time $n$ is that at least one of
the four branches $C_i$ contains a path of length at least $n$
\emph{and} no stable site.
Since a protected site is stable, to complete the proof,
it suffices to show that there are $\alpha>0$ and $K < \infty$
such that
\begin{equation} \label{expon}
P_p(|C_i| \geq n \text{ and $C_i$
contains no protected site}) \leq K \, e^{- \alpha n}.
\end{equation}
To prove (\ref{expon}), we partition ${\mathbb Z}^2$
into disjoint plaquettes and denote by $S$ the collection
of these plaquettes.
We do an algorithmic construction of $C_i$
(as in, e.g., \cite{fn}), where the order of
checking the state of sites is such that when the
first site in a square from $S$ is checked and found
to be closed, then the other three sites in that
plaquette are checked next.
Then standard arguments show that the probability
in~(\ref{expon}) is bounded above by
$K \, [1 - (1 - p)^4]^{(n /4)}$. \fbox{}

\begin{remark} \label{remark1}
We note that one can improve Theorem~\ref{convergence},
namely prove exponential convergence \emph{uniformly} in
$p \in [0,1]$ (i.e., it is possible to get a constant
$c_0>0$ independent of $p$).
This is done by using the proof given above for values
of $p$ smaller than some $p_0>p_c^*$, together with the fact
that for $p \geq p_0$ the size of the closed cluster of
the origin at time $0$ has an exponential tail~\cite{grimmett}.
(We have chosen to give the argument in the proof simply
because it has the advantage of being valid for all values
of $p$.)
\end{remark}

Using arguments analogous to those in the proof of
Theorem~\ref{convergence}, one can get exponential decay
of correlations for $\bar P_p$, which was proved, in a
somewhat different way, in~\cite{fis}.
This result is important in the context of the present
paper because it suggests (see, for example,~\cite{cardy2})
that $\bar P_p$ is in the universality class of independent
percolation, as we will show in the next section.
We include here a proof of the result (see~\cite{fis} for
the original proof) in order to show how the same mechanism
is responsible for the exponential decay of correlations of
the limiting measure and the exponential convergence to the
final configuration (Theorem~\ref{convergence}).
As it will be clear from the proofs of the main results,
that same mechanism is also responsible for the fact that
$\bar P_p$ is in the universality class of independent
percolation.
Such mechanism explains, in this particular case, the
model's ``universal" behavior and its relation with the
exponential decay of correlations.

For $x \in {\mathbb Z}^2$, let $\text{d}(o,x)$ be one plus
the number of sites between $o$ and $x$ along a shortest
${\mathbb Z}^2$-path from $o$ to $x$, and let
$B_x(r) = \{ y \in {\mathbb Z}^2 : \text{d}(y,x) < r \}$.

\begin{theorem} \emph{\cite{fis}} \label{fis}
$\bar P_p$ has exponentially decaying correlations:
\begin{equation} \label{exp-corr}
|E_p[\bar\omega(o) \bar\omega(x)] - E_p[\bar\omega(o)] E_p[\bar\omega(x)]|
\leq R \,\exp{[-c'_0 \, \text{\emph{d}}(o,x)]},
\end{equation}
where $R < \infty$ and $c_0' > 0$.
\end{theorem}

\noindent \emph{Proof.} Denote by $A_x(n)$ the event that
$\bar\omega(x)$ is determined only by the configuration
$\omega_{B_x(n)}$ inside $B_x(n)$ and by $A^c_x(n)$ its complement.
The proof rests on the observation that if $\text{d}(o,x)>2n$,
then conditioned on $A_o(n)$ and $A_x(n)$, the random variables
$\bar\omega(o)$ and $\bar\omega(x)$ are independent.

Before proceeding with the proof, we notice that a necessary
condition for $A^c_o(n)$ to occur is that the origin be closed
at time $0$ and that there be at least one branch $C_i$ of the
cluster of the origin at time $0$ that reaches the boundary of
$B_o(n)$ and has no stable site inside $B_o(n)$.
This event is analogous to the one considered at the end of the
proof of Theorem~\ref{convergence}.
Then, arguments analogous to those used there to get~(\ref{expon})
give the bound
\begin{equation} \label{bound4}
P_p[A^c_o(n)] \leq \exp{(- \alpha' \, n)},
\end{equation}
for some $\alpha'>0$.

Take $N$ such that $P_p[A^c_o(N)] < 1/2$ and consider the set of sites
$ \{x \in {\mathbb Z}^2 : \text{d}(o,x) \geq 3N \} = {\mathbb Z}^2
\setminus B(3N)$.
For a site in ${\mathbb Z}^2 \setminus B(3N)$, we take
$n = \lceil \text{d}(o,x)/3 \rceil$ and write, thanks to the
observation above,
\begin{eqnarray}
E_p[\bar\omega(o) \, \bar\omega(x)] & = &
E_p[\bar\omega(o) \, \bar\omega(x) \, | \, A_o(n) \cap A_x(n)] \, \{1 - P_p[A^c_o(n)] \}^2 \nonumber \\
 & + & E_p[\bar\omega(o) \, \bar\omega(x) \, | \, A^c_o(n) \cup A^c_x(n)] \, \{2 - P_p[A^c_o(n)]\} \, P_p[A^c_o(n)] \\
 & = & E_p[\bar\omega(o) \, | \, A_o(n) \cap A_x(n)] \,
E_p[\bar\omega(x) \, | \, A_o(n) \cap A_x(n)] \, \{1 - P_p[A^c_o(n)]\}^2  \nonumber \\
 & + & \, E_p[\bar\omega(o) \, \bar\omega(x) \, | \, A^c_o(n) \cup A^c_o(n)] \,
\{2 - P_p[A^c_o(n)]\} \, P_p[A^c_o(n)], \label{exp-decay1}
\end{eqnarray}
where we have used
\begin{equation}
P_p[A_o(n) \cap A_x(n)] = P_p[A_o(n)] P_p[A_x(n)] = \{1 - P_p[A^c_o(n)]\}^2
\end{equation}
and
\begin{equation}
P_p[A^c_o(n) \cup A^c_x(n)] = 1 - P_p[A_o(n) \cap A_x(n)] = \{2 - P_p[A^c_o(n)]\} \, P_p[A^c_o(n)] ,
\end{equation}
which follow from the observation that $A_o(n)$ and $A_x(n)$
are independent events because $\text{d}(o,x)>2n$.

We now write
\begin{equation}
E_p[\bar\omega(o) \, | \, A_o(n) \cap A_x(n)] = 
\frac{ E_p[\bar\omega(o)] - E_p[\bar\omega(o) \, | \, A^c_o(n) \cup A^c_x(n)] \,
\{2 - P_p[A^c_o(n)]\} \, P_p[A^c_o(n)] }{ \{1 - P_p[A^c_o(n)]\}^2 }
\end{equation}
and the same for $E_p[\bar\omega(x) \, | \, A_o(n) \cap A_x(n)]$,
and plug the two expressions in~(\ref{exp-decay1}) to get
\begin{equation}
E_p[\bar\omega(o) \bar\omega(x)]
= \frac{1}{\{1 - P_p[A^c_o(n)]\}^2} \, E_p[\bar\omega(o)] \, E_p[\bar\omega(x)]
+ R_1 \, P_p[A^c_o(n)] + R_2 \, P_p[A^c_o(n)]^2, \label{exp-decay2}
\end{equation}
for some constants $R_1$ and $R_2$.
From~(\ref{bound4}) and~(\ref{exp-decay2}), we immediately see that
\begin{equation}
|E_p(\bar\omega(0) \, \bar\omega(x)) - E_p(\bar\omega(o)) \, E_p(\bar\omega(x))|
\leq R \, P_p[A^c_o(n)] \leq R \, e^{-c \, n},
\end{equation}
for some $R < \infty$ and $c>0$.

For the sites in ${\mathbb Z}^2 \setminus B(3N)$, the proof
is concluded by taking $c_0' = c / 3$.
For the sites in $B(3N)$, we just have to choose a constant
$R$ large enough so that $R \, \exp{(-3Nc_0')} \geq 1$. \fbox{}

\begin{remark} \label{remark2}
%
We note that one can get exponential decay of correlations
\emph{uniformly} in $p \in [0,1]$ (as in~\cite{fis}) using the
fact that for $p > p_c^*$ the size of the closed cluster
of the origin at time $0$ has an exponential tail~\cite{grimmett}
(see Remark~\ref{remark1}).
\end{remark}

We conclude this section with Proposition~\ref{critical}, which
identifies the critical density of our bootstrap percolation model
on ${\mathbb Z}^2_{cp}$ with $p_c^*$, showing that the bootstrapping
rule employed here does not shift the critical point.
This motivates the next section, where we analyze the continuum scaling
limit of crossing probabilities and some critical exponents of the
bootstrapped model on ${\mathbb Z}^2_{cp}$ when the initial density
of $1$'s is $p_c^*$.

\begin{proposition} \label{critical}
The following results hold for $\bar\omega$.
\begin{itemize}
\item[1.]
\begin{enumerate}
\item[(i)] Closed sites do not percolate if $p>p_c^*$ and
percolate if $p<p_c^*$.
\item[(ii)] If $p=p_c^*$, closed sites do not percolate
and the mean cluster size for the closed component is infinite.
\end{enumerate}
\item[2.]
\begin{enumerate}
\item[(i)] Open sites do not $*$-percolate if $p<p_c^*$
and $*$-percolate if $p>p_c^*$.
\item[(ii)] If $p=p_c^*$, open sites do not $*$-percolate
and the mean $*$-cluster size for the open component
is infinite.
\end{enumerate}
\end{itemize}
\end{proposition}

\noindent {\it Proof}. Let us begin with the proofs of {\it 1.(i)}
and {\it 2.(i)}, which are elementary.
If $p>p_c^*$, closed sites do not percolate in $\omega$, that is before
bootstrapping the open sites, and therefore cannot possibly percolate
in $\bar\omega$, after bootstrapping the open sites.
If $p<p_c^*$, on the contrary, closed sites do percolate in $\omega$,
and since any doubly-infinite closed ${\mathbb Z}^2$-path (i.e., a
closed ${\mathbb Z}^2$-path that can be split in two disjoint infinite
paths) contained in $\omega$ is stable and therefore it is also contained
in $\bar\omega$, this implies that closed sites percolate in $\bar\omega$
and concludes the proof of {\it 1.(i)}.

To prove {\it 2.(i)}, it suffices to notice that for $p<p_c^*$, closed sites
percolate in $\omega$ and the origin is surrounded by infinitely many
${\mathbb Z}^2$-loops of closed sites.
Such closed loops are stable and therefore still exist in $\bar\omega$
and prevent open sites from $*$-percolating.
On the other hand, if $p>p_c^*$ open sites $*$-percolate already in $\omega$,
which concludes the proof of {\it 2.(i)}.

{\it 1.(ii)} and {\it 2.(ii)} can be proved together using a theorem
of Russo~\cite{russo}.
At $p=p_c^*$, in $\omega$ the origin is surrounded by infinitely many
${\mathbb Z}^2$-loops of closed sites and infinitely many $*$-loops of
open sites.
Both types of loops are stable and therefore in $\bar\omega$ there is
no percolation of closed sites, nor $*$-percolation of open sites.
By an application of a theorem of Russo~\cite{russo}, this implies that
both the mean cluster size of the closed component and the mean $*$-cluster
size of the open component diverge. \fbox{} \\

\section{Main Results} \label{main-results}

In this section we present the main results of this paper;
the proofs will be given in Section~\ref{proofs}.
The results presented in this section hold for all the measures
that are intermediate between the initial measure $P_p$ and the limiting
one $\bar P_p$.
These form a one parameter family $\{ P_{p,n} \}_{n \in {\mathbb N}}$
of measures, parametrized by time $n = 1, 2, 3, \ldots$, and are
increasingly different from $P_p$ as $n$ becomes larger.

\subsection{The Continuum Scaling Limit of Crossing Probabilities}
\label{universality-scaling}

We take a ``mesh'' $\delta$ and consider the ``scaling
limit'' of crossing probabilities for the percolation model
$\bar \omega$ on $\delta {\mathbb Z}^2$ as $\delta \to 0$,
focusing for simplicity on the probability of an open
$*$-crossing of a rectangle aligned with the coordinate axes.
A similar approach would work for any domain with a ``regular''
boundary, but it would imply dealing with more complex deformations
of the boundary than that needed for proving the result for a
rectangle.

Consider a finite rectangle
${\cal R} = {\cal R}(b,h) \equiv (-b/2,b/2) \times (-h/2,h/2)
\subset {\mathbb R}^2$ centered at the origin of ${\mathbb Z}^2$,
with sides of lengths $b$ and $h$ and aspect ratio $\rho = b/h$.
We say that there is an open vertical $*$-crossing of ${\cal R}$
in $\omega$ (resp., $\bar\omega$) if ${\cal R} \cap
\delta{\mathbb Z}^2$ contains a $*$-path of open sites from
$\omega$ (resp., $\bar\omega$) joining the top and bottom sides
of the rectangle ${\cal R}$, and call $\phi^*_{\delta}(b,h;n)$ the
probability of such an open crossing at time $n$.

More precisely, there is a vertical open $*$-crossing at time $n$
if there is a $*$-path $(x_0, x_1, \ldots, x_m, x_{m+1})$ in
${\mathbb Z}^2$ such that $\omega_n(x_j)=1$ for all $j$,
$\delta x_0, \delta x_1, \ldots, \delta x_m, \delta x_{m+1}$
are all in $\cal R$, and the line segments $\overline{\delta x_0,\delta x_1}$
and $\overline{\delta x_m,\delta x_{m+1}}$ touch respectively
the top side $[-b/2,b/2] \times \{ h/2 \}$ and the bottom side
$[-b/2,b/2] \times \{ -h/2 \}$ of $\cal R$.

It is believed that the scaling limit of crossing probabilities
for \emph{independent} percolation exists and is given by Cardy's
formula (see~\cite{cardy1,cardy2}); this has however been rigorously
proved only for critical site percolation on the triangular
lattice~\cite{smirnov2}.
We will \emph{assume} that $\lim_{\delta \to 0} \phi^*_{\delta}(b,h;0)
= F(\rho)$, where $F$ is a continuous function of its argument.

\begin{theorem} \label{cross-prob}
Suppose that the scaling limit of the crossing probability
of a rectangle $\cal R$ exists for independent critical site
percolation on ${\mathbb Z}^2_{cp}$ and is given by a continuous
function $F$ of $\rho$.
Then, the corresponding crossing probability in the bootstrapped
model $\bar\omega$ with $p = p_c^*$ has the same scaling limit.
\end{theorem}

\subsection{Critical Exponents} \label{critexp}

We will consider four percolation critical exponents, namely
the exponents $\beta$ (related to the percolation probability), $\nu$
(related to the correlation length), $\eta$ (related to the connectivity
function) and $\gamma$ (related to the mean cluster size).
The existence of these exponents has been recently proved~\cite{lsw,sw},
and their predicted values confirmed rigorously, for the case of independent
site percolation on the triangular lattice.
Such exponents are believed to be universal for independent percolation in
the sense that their value should depend only on the number of dimensions and not
on the structure of the lattice or on the nature of the percolation model (e.g.,
whether it is site or bond percolation); that type of universality has not yet
been proved.

Consider an independent percolation model with distribution $P_p$ on a
two-dimensional lattice $\mathbb L$ such that $0<p_c<1$.
Let $C_o$ be the open cluster containing the origin and $|C_o|$ its cardinality,
then $\theta(p) = P_p(|C_o| = \infty)$ is the \emph{percolation probability}.
Arguments from theoretical physics suggest that $\theta(p)$ behaves roughly like
$(p-p_c)^{\beta}$ as $p$ approaches $p_c$ from above.

It is also believed that the \emph{connectivity function}
\begin{equation} \label{connect-function}
\tau_p(x) = P_p(\text{the origin and } x \text{ belong to the same cluster})
\end{equation}
behaves, for the Euclidean length $||x||$ large, like $||x||^{-\eta}$ if $p=p_c$,
and like $\exp{(-||x||/\xi(p))}$ if $0<p<p_c$, for some $\xi(p)$ satisfying
$\xi(p) \to \infty$ as $p \uparrow p_c$.
The \emph{correlation length} $\xi(p)$ is defined by
\begin{equation} \label{correlation-length}
\xi(p)^{-1} = \lim_{||x|| \to \infty} \left\{ - \frac{1}{||x||} \log \tau_p(x) \right\}.
\end{equation}
$\xi(p)$ is supposed to behave like $(p_c-p)^{-\nu}$ as $p \uparrow p_c$.
The \emph{mean cluster size} $\chi(p) = E_p |C_o|$ is also believed to diverge
with a power law behavior $(p_c-p)^{-\gamma}$ as $p \uparrow p_c$.

It is not clear how strong one may expect such asymptotic relations to be
(for more details about critical exponents and \emph{scaling theory} in
percolation, see~\cite{grimmett} and references therein); for this reason
the logarithmic relation is usually employed.
This means that the previous conjectures are usually stated in the following form:
\begin{eqnarray}
\lim_{p \downarrow p_c} \frac{\log \theta(p)}{\log (p-p_c)} = \beta, \\
\lim_{||x|| \to \infty} \frac{\log \tau_{p_c}(x)}{\log ||x||} = - \eta, \\
\lim_{p \uparrow p_c} \frac{\log \xi(p)}{\log (p_c - p)} = - \nu, \\
\lim_{p \uparrow p_c} \frac{\log \chi(p)}{\log (p_c - p)} = - \gamma.
\end{eqnarray}

In the rest of the paper, $\theta(p)$, $\tau_p(x)$, $\xi(p)$ and $\chi(p)$
will indicate the percolation probability, connectivity function, correlation
length and mean cluster size for independent site percolation on ${\mathbb Z}^2_{cp}$.
For $n \in [1, \infty]$, let $\theta(p,n)$, $\tau_{p,n}(x)$, $\xi(p,n)$ and
$\chi(p,n)$ be respectively the percolation probability, connectivity function,
correlation length and mean cluster size on ${\mathbb Z}^2_{cp}$ for the
bootstrapped model at time $n$, with $n=\infty$ corresponding to the fully
bootstrapped configuration $\bar\omega$.
The main theorem of this section is the following.

\begin{theorem} \label{mainthm-critexp}
There exist constants $0 < c_1, c_2 < \infty$ such that,
$\forall n \in [1, \infty]$,
\begin{eqnarray}
\theta(p) \leq \theta(p,n) \leq c_1 \, \theta(p), \,\,\,
\text{ for } p \in (p^*_c,1], \label{theta} \\
\tau_p(x) \leq \tau_{p,n}(x) \leq p^{-c_2} \, \tau_p(x), \,\,\,
\text{ for } p \in (0, p^*_c], \label{tau} \\
\xi(p,n) = \xi(p), \,\,\, \text{ for } p \in (0, p^*_c]. \label{xi}
\end{eqnarray}
\end{theorem}

The next corollary is an immediate consequence of Theorem~\ref{mainthm-critexp}
and its main application; it says that the bootstrapped percolation model
(in fact, all models corresponding to $n$ enhancements by bootstrapping,
with $n \in [1,\infty]$) has the same critical exponents $\beta$, $\eta$,
$\nu$ and $\gamma$ as ordinary independent percolation.

\begin{corollary} \label{cor-critexp}
Suppose that the critical exponents  $\beta$, $\eta$, $\nu$ and $\gamma$
exist for independent site percolation on ${\mathbb Z}^2_{cp}$, then they
also exist for the bootstrapped model and have for the latter the same
numerical values as for the original model.
\end{corollary}

\section{Proofs of the Main Results} \label{proofs}

In this section we prove the main results of this paper,
presented in Section~\ref{main-results}.

\subsection{Crossing Probabilities -- Proof of Theorem~\ref{cross-prob}}

To prove the theorem we need to compare the probability of an open
vertical $*$-crossing of $\cal R$ in $\bar\omega$ with the probability of
the same event in $\omega$.
In order to do that, we will use the natural coupling that exists
between $\omega$ and $\bar\omega$ via bootstrapping.
First of all notice that, if an open vertical $*$-crossing of
$\cal R$ is present in $\omega$, it is also present in $\omega_n$,
for all $n$, since open sites are stable.
Therefore,
\begin{equation} \label{lower-bound}
\lim_{\delta \to 0} \phi^*_{\delta}(b,h;n) \geq
\lim_{\delta \to 0} \phi^*_{\delta} (b,h;0) = F(\rho).
\end{equation}
(\ref{lower-bound}) holds for all values of $n$, including
$n = \infty$, so if we call $\bar \phi^*_{\delta}(b,h)$
the probability of an open vertical $*$-crossing of $\cal R$
from $\bar \omega$, we can write
\begin{equation} \label{lower-bound1}
\lim_{\delta \to 0} \bar \phi^*_{\delta}(b,h) \geq
\lim_{\delta \to 0} \phi^*_{\delta} (b,h;0) = F(\rho).
\end{equation}

On the other hand, if an open vertical $*$-crossing of $\cal R$
is not present in $\omega$, this implies the existence of a
\emph{closed horizontal ${\mathbb Z}^2$-crossing} of $\cal R$.
For $\delta$ small such a crossing must involve many sites,
and the probability of finding ``near'' its endpoints two sites
$x$ and $y$, belonging to the crossing, attached through closed
${\mathbb Z}^2$-paths to two stable closed sites $x'$ and $y'$
should be close to one.
If such stable sites are found, Lemma~\ref{lemma-protected}
assures that at least the portion of the closed horizontal
crossing from $x$ to $y$ is still present in $\bar\omega$.
This suggests that, conditioned on having in $\omega$
a closed horizontal ${\mathbb Z}^2$-crossing of a slightly
bigger (in the horizontal direction) rectangle, with high
probability, in $\bar\omega$ there will be a closed horizontal
${\mathbb Z}^2$-crossing of $\cal R$ blocking any open
vertical $*$-crossing.
It is then enough to prove that this probability goes to one
as $\delta \to 0$.

We will now make this more precise, adapting the
proof of Theorem~1 of~\cite{cns1}.
Consider the rectangle ${\cal R}'= {\cal R}(b',h)$
with $b'$ slightly larger than $b$ and aspect ratio
$\rho'=b'/h$.
It follows from our assumptions that
\begin{equation} \label{lim1}
\phi^*(b',h;0) \equiv
\lim_{\delta \to 0} \phi^*_{\delta} (b',h;0) = F(\rho')
\end{equation}
and
\begin{equation} \label{lim2}
\lim_{b' \to b} \phi^*(b',h;0) =
\lim_{\rho' \to \rho} F(\rho') = F(\rho).
\end{equation}

If we now call $\phi_{\delta}(b',h+\delta;0)$ the
probability of a \emph{closed horizontal ${\mathbb Z}^2$-crossing}
of ${\cal R}(b',h+\delta)$ from $\omega$, and
$\phi_{\delta}(b,h+\delta;n)$ that of a
\emph{closed horizontal ${\mathbb Z}^2$-crossing} of
${\cal R}(b,h+\delta)$ from $\omega_n$, the observation
that a closed ${\mathbb Z}^2$-crossing can only be
``eaten'' from its endpoints yields
\begin{equation} \label{bound}
\phi_{\delta} (b,h+\delta;n) \geq \phi_{\delta} (b',h+\delta;0),
\end{equation}
as long as $b'>b$, and $n$ is not too large (depending on
$b'-b$ and $\delta$).

Since a closed horizontal ${\mathbb Z}^2$-crossing
of ${\cal R}(b,h+\delta)$ blocks any open vertical
$*$-crossing of ${\cal R}(b,h)$ and vice versa,
(\ref{bound}) yields
\begin{equation} \label{upper-bound}
\phi^*_{\delta} (b,h;n)
= [ 1 - \phi_{\delta} (b,h+\delta;n) ]
\leq [ 1 - \phi_{\delta} (b',h+\delta;0) ]
= \phi^*_{\delta} (b',h;0).
\end{equation}
Keeping $n$ fixed, we can let first $\delta$ go to zero
and then $b'$ go to $b$, thus obtaining from~(\ref{upper-bound})
a bound that, combined with~(\ref{lower-bound}), gives the
desired result, at least for values of $n$ that are not too large.

To complete the proof, we will extend~(\ref{bound})
to all values of $n$, including $n = \infty$, at the cost
of a correction that goes to zero as $\delta \to 0$.
In order to do this, we will use Lemma~\ref{stable-path}
to show that if there is a closed horizontal crossing by
$(y_0, \dots, y_k)$ of ${\cal R}(b',h+\delta)$ at time $0$,
with high probability it does not ``shrink'' too much due
to the effect of the dynamics, so that at all later times,
including $n = \infty$, there is a closed horizontal
crossing of ${\cal R}(b,h+\delta)$ by $(y_{k_1}, \dots, y_{k_2})$.
This is achieved by looking at the partial clusters
containing the portions of $(y_0, \dots, y_k)$
contained in ${\cal R}(b',h+\delta) \setminus {\cal R}(b,h+\delta)$
and searching for protected sites.

Noting that each of the partial paths
$(y_0, \dots, y_{k_1})$ and $(y_{k_2},\dots,y_k)$
contains of the order of $(b'-b)/\delta$ sites,
we see that Lemma~\ref{stable-path} implies that
it suffices to show that there exist $\alpha>0$
and $K < \infty$ such that for any deterministic
$(x,x')$,
\begin{equation} \label{expo}
P_{p_c^*}(|C_{(x,x')}| \geq \ell \text{ and $C_{(x,x')}$
contains no protected site}) \leq K \, e^{- \alpha \ell}.
\end{equation}
To prove (\ref{expo}), we proceed as in the proof
of Theorem~\ref{convergence}, that is, we
partition ${\mathbb Z}^2$ into disjoint plaquettes
and denote by $S$ the collection of these plaquettes.
We then do an algorithmic construction of $C_{(x,x')}$
where the order of checking the state
of sites is such that when the first site in a plaquette
from $S$ is checked and found to be closed, then the
other three sites in that plaquette are checked next.
Again, standard arguments show  that the probability
in~(\ref{expo}) is bounded above by
$K \, [1 - (1 - p_c^*)^4]^{(\ell /4)}$. \fbox \\

\begin{remark}
As already remarked, the proof of Theorem~\ref{cross-prob}
shows that the result is valid for all the intermediate
measures $P_{p_c^*,n}$.
\end{remark}

\subsection{Critical exponents} \label{critexp-proofs}

\subsubsection{Proof of Theorem~\ref{mainthm-critexp}}

For two subsets $C$ and $D$ of ${\mathbb Z}^2$, we denote
by $\{ C \longleftrightarrow D \}$ the event that some site
in $C$ is connected to some site in $D$ by an open $*$-path,
and by $\{ C \longleftrightarrow \infty \}$ the event that
some site in $C$ belongs to an infinite open $*$-path.

The lower bound for $\theta(p,n)$ is obvious.
For the upper bound, we let ${\cal N}^*_x$ be the set of
$*$-neighbors of $x$ and rely on the following observation.
If no site in ${\cal N}^*_o$ belongs to an infinite open
$*$-path at time $0$, then the origin must be surrounded
by a closed ${\mathbb Z}^2$-loop $\lambda$.
It then follows, by Lemma~\ref{lemma-protected}, that
each site in $\lambda$ is stable.
Therefore, the origin will not be connected to infinity
by an open $*$-path at any later time.
Thus,
\begin{equation}
\theta(p,n) \leq P_p({\cal N}^*_o \longleftrightarrow \infty).
\end{equation}

Since $\{ o \longleftrightarrow \infty \}$ can be written as
$\{ \omega(o)=1 \} \cap \{ {\cal N}^*_o \longleftrightarrow \infty \}$,
and $\{ \omega(o)=1 \}$ and $\{ {\cal N}^*_o \longleftrightarrow \infty \}$
are independent at time $0$,
\begin{equation}
P_p(o \longleftrightarrow \infty) = p \, P_p({\cal N}^*_o \longleftrightarrow \infty).
\end{equation}
From this  we get
\begin{equation}
\theta(p,n) \leq p^{-1} \, \theta(p) \leq \frac{1}{p_c^*} \, \theta(p),
\end{equation}
as required.

The lower bound for $\tau_{p,n}(x)$ is again obvious.
To obtain the upper bound, we first note that for $||x||$ bounded,
the inequality is trivial by choosing $c_2$ big enough so that the
right-hand side of~(\ref{tau}) exceeds $1$.
Next, for $||x||$ large enough, we notice that, unless
$\{ {\cal N}^*_o \longleftrightarrow {\cal N}^*_x \}$ at time $0$,
the origin and $x$ must be separated by a closed ${\mathbb Z}^2$-loop
surrounding one of them or by a doubly-infinite closed ${\mathbb Z}^2$-path,
and therefore it cannot be the case that $\{ 0 \longleftrightarrow x \}$
at any later time.
Thus,
\begin{equation}
\tau_{p,n}(x) \leq P_p({\cal N}^*_o \longleftrightarrow {\cal N}^*_x).
\end{equation}

Since $\{o \longleftrightarrow x\}$ can be written as
$\{ \omega(o)=\omega(x)=1 \} \cap
\{ {\cal N}^*_o \longleftrightarrow {\cal N}^*_x) \}$,
and $\{ \omega(o)=\omega(x)=1 \}$ and
$\{ {\cal N}^*_o \longleftrightarrow {\cal N}^*_x) \}$
are independent at time $0$,
\begin{equation}
P_p(o \longleftrightarrow x) = p^2 \, P_p({\cal N}^*_o \longleftrightarrow {\cal N}^*_x).
\end{equation}
From this  we get
\begin{equation}
\tau_{p,n}(x) \leq p^{-2} \, \tau_p(x),
\end{equation}
as required.

Eq.~(\ref{xi}) is an immediate consequence of~(\ref{tau}) and the
definition of $\xi(p)$; it is enough to observe that
\begin{equation}
\lim_{||x|| \to \infty} \left\{ - \frac{1}{||x||} \left[\log \tau_p(x) - c_2 \log p \right] \right\}
= \xi(p)^{-1}. \,\,\,\,\, \fbox{}
\end{equation}

\subsubsection{Proof Corollary~\ref{cor-critexp}}

It follows from~(\ref{theta})~and~(\ref{tau}) that, for $p \in (p_c^*,1]$ and $||x||>1$,
\begin{eqnarray}
- \frac{\log \theta(p)}{\log (p - p_c^*)} \leq - \frac{\log \theta(p,n)}{\log (p - p_c^*)}
\leq - \frac{\log \theta(p) + \log c_1}{\log (p - p_c^*)}, \label{bound1} \\
\frac{\log \tau_{p_c^*}(x)}{\log ||x||} \leq \frac{\log \tau_{p_c^*,n}(x)}{\log ||x||}
\leq \frac{\log \tau_{p_c^*}(x) - c_2 \log p_c^*}{\log ||x||}. \label{bound2}
\end{eqnarray}
For $p \in (0,p_c^*)$, observing that
$\chi(p) = E_p \sum_{x \in {\mathbb Z}^2} I(o \longleftrightarrow x)
= \sum_{x \in {\mathbb Z}^2} \tau_p(x)$ (where $I(\cdot)$ is the indicator function),
(\ref{tau}) yields $\chi(p) \leq \chi(p,n) \leq p^{- \, c_2} \chi(p)$, and therefore
\begin{equation}
- \frac{\log \chi(p)}{\log (p - p_c^*)} \leq - \frac{\log \chi(p,n)}{\log (p - p_c^*)}
\leq - \frac{\log \chi(p) - c_2 \log p}{\log (p - p_c^*)}. \label{bound3}
\end{equation}

Using (\ref{bound1}), (\ref{bound2}) and (\ref{bound3}), together with~(\ref{xi})
and the definitions of the critical exponents, and taking the appropriate limits
gives the desired results. \fbox{}

\bigskip
\bigskip

\noindent {\bf Acknowledgements.} The author thanks M.~Isopi
for starting his interest in the model considered here and
for pointing at ref.~\cite{fis}, and an anonymous referee
for helpful comments on the presentation of the paper.
Interesting and helpful discussions with L.~R.~Fontes,
F.~den~Hollander, M.~Isopi, A.~Sakai, V.~Sidoravicius
and A.~C.~D.~van~Enter are also acknowledged.

\end{document}